\documentclass[11pt,a4paper]{article}
\usepackage{oldlfont}

\usepackage{amsmath}
\usepackage{amsfonts}
\usepackage{amssymb}
\usepackage{array,longtable}
\usepackage{amscd}

\begin{document}

\begin{center}
\rule{0pt}{4cm}
    {\huge {\bf Gantmakher--Krein theorem \\
    for 2-totally nonnegative operators \\[0.3cm]
    in ideal spaces}} \\[1cm]
\end{center}
\begin{center}
    {\Large Olga Y. Kushel, Petr P. Zabreiko \\
Department of mechanics and mathematics, \\
Belorussian State University, \\ Nezavisimosti sq., 4, 220050,
Minsk, Belarus, \\[0.15cm]
 e-mail: kushel@mail.ru;
zabreiko@bsu.by}
\end{center}
\medskip

\begin{center}
{\Large August 2008}
\end{center}

\bigskip

\begin{center}
{\bf Abstract.}
\end{center}

The tensor and exterior squares of a completely continuous
non-negative linear operator $A$ acting in the ideal space
$X(\Omega)$ are studied. The theorem representing the point
spectrum (except, probably, zero) of the tensor square $A \otimes
A$ in the terms of the spectrum of the initial operator $A$ is
proved. The existence of the second (according to the module)
positive eigenvalue $\lambda_2$, or a pair of complex adjoint
eigenvalues of a completely continuous non-negative operator $A$
 is proved under the additional condition, that its exterior
square $A\wedge A$ is also nonnegative.

\medskip


\medskip


\medskip

\section{Introduction} This paper presents
the results of our study of the spectrum of 2-totally-nonnegative
operators in ideal spaces. The theory of total positivity is
mainly based on the results of F.R. Gantmakher and M.G. Krein (see
[1]), concerning the properties of the spectrum of {\it
$k$-totally primitive} matrices (i.e. matrices, which are
nonnegative and primitive along with their $j$-th compound
matrices $(1 < j \le k)$ up to the order $k$). In the most
important case $k = n$ $k$-totally primitive matrices are called
{\it oscillatory}.
 In monograph [1] the following
statement was proved: {\it if the matrix $\mathbf{A}$ of a linear
operator $A$ in the space ${\mathbb{R}}^{n}$ is $k$-totally
primitive, then the operator $A$ has $k$ positive simple
eigenvalues $0 < \lambda_k < \ldots < \lambda_2 < \lambda_1$, with
a positive eigenvector $e_1$ corresponding to the maximal
eigenvalue $\lambda_1$, and an eigenvector $e_j$, which has
exactly $j - 1$ changes of sign, corresponding to $j$-th
eigenvalue $\lambda_j$} (see [1], p. 310, theorem 9).

The study of linear integral operators with analogous properties
predates the study of oscillatory matrices. The main results
concerning this problem were received by O.D. Kellog (see [2]). He
proved the theorem about spectral properties of continuous
symmetric totally nonnegative kernels. Later this theorem was
generalized by F.R. Gantmakher for the non-symmetric case. This
result one can find in monograph [1] in the following form: {\it
let $k(t,s) \in C[0,1]^2$ satisfy the following conditions:

(a) for any $0< t_0 < t_1 < \ldots < t_n < 1$ and $0< s_0 < s_1 <
\ldots < s_n < 1$ $ \ \ n = 0, \ 1, \ \ldots$ the inequality
$$k\begin{pmatrix}
  t_0 & t_1 & \ldots & t_n \\
  s_0 & s_1 & \ldots & s_n
\end{pmatrix} \geq 0$$
is true;

 (b) for any $0< t_0 < t_1 < \ldots < t_n <
1$ $ \ \ n = 0, \ 1, \ \ldots$ the inequality
 $$k\begin{pmatrix} t_0 & t_1 & \ldots & t_n \\
  t_0 & t_1 & \ldots & t_n
\end{pmatrix} > 0$$ is true.

Then all the eigenvalues of the linear integral equation
$$\int_0^1k(t,s)x(s)ds = \lambda x(t)$$ are positive and simple:
$$0 < \ldots < \lambda_n < \ldots < \lambda_2 < \lambda_1,$$
 with
a strictly positive on $(0,1)$ eigenfunction $e_1(t)$
corresponding to the maximal eigenvalue $\lambda_1$, and an
eigenfunction $e_n(t)$, which has exactly $n - 1$ changes of sign
and no other zeroes on $(0,1)$, corresponding to the $n$-th
eigenvalue $\lambda_n$} (see [1], p. 211).

The proof of this statement one can find also in [3], where the
history of the theory of totally positive matrices and kernels is
presented in detail. Unlike monograph [1], in which the basic
statements of the theory are given in the form most suitable for
the study of  small oscillations of mechanical systems, in [3]
definitions and theorems about the properties of totally positive
kernels are given in the pure form.

 In paper [4] by S.P. Eveson the  result mentioned  was spread onto
a wider class of kernels. The existence of $k$ positive
eigenvalues was proved under some additional assumptions for the
case of a compact linear integral operator, acting in $L_2[0,1]$,
which kernel is totally positive of order $k$. A substantial
contribution into the development of the theory of totally
positive and sign-symmetric kernels was made by S. Karlin (see
[5]).

 Once a great number of papers are devoted to
the theory of totally positive matrices and kernels,  in the case
of abstract (not necessarily integral) compact linear operators
the situation is absolutely different. Here we  can mention only a small
number of papers. In paper [6] oscillatory operators in $C[a,b]$
were studied by the method of the passage to the limit from
finite-dimensional approximants. In paper [7] another method of
generalization was suggested. But this method was realized also
only for the space $C[a,b]$. Many results, related to the
applications of the theory of osscillation to differential
operators, were included into monograph [8] by Yu.V. Pokornyi and
his group.

In paper [9] we studied $2$-totally indecomposable operators (i.e. indecomposable
operators that are
nonnegative with respect to some cone $K$, and such that their  exterior squares are also nonnegative and
indecomposable) in the spaces $L_p(\Omega) \ (1 \le p \le \infty)$
and $C(\Omega)$. We proved the existence of the second (according
to the module) eigenvalue $\lambda_2$ of a completely continuous
non-negative operator $A$ under the condition that its exterior
square $A\wedge A$ is also non-negative. The simplicity of the
first and second eigenvalues was proved and the interrelation
between the indices of imprimitivity of $A$ and $A\wedge A$ was
examined for the case, when the operators $A$ and $A\wedge A$ are
indecomposable. The difference (according to the module) of
$\lambda_1$ and $\lambda_2$ from each other and from other
eigenvalues was proved for the case, when $A$ and $A\wedge A$ are
primitive.

 In the present paper we are going
to generalize the results, received in [9], for $2$-totally
nonnegative operators in some ideal spaces.  As the authors
believe, the natural method of the examination of such operators
is a crossway from studying an operator $A$, acting in an ideal
space $X$ to the study of the operators $A \otimes A$ and $A \wedge
A$, acting in spaces with mixed norms. Let us turn now to a more
detailed outline of the paper. In Section 2 we briefly consider
the basic properties of ideal spaces. Tensor and exterior squares of
ideal spaces are described in Section 3. The connection between
the topological exterior square of the ideal space $X$ and the
space with mixed norm is also clarified in Section 3. In Section 4
we work with the tensor and exterior square of a linear operator
$A: X \rightarrow X$. These operators act  in the tensor and
exterior square of the initial space $X$ respectively. Generally, in Sections 3
and 4 we develop the mathematical tools that will enable us to
define the class of abstract totally nonnegative operators and to
generalize the results of Gantmakher and Krein. Sections 5 and 6
present a number of  results on  the description of the spectrum and the
parts of the spectrum of the tensor square and the exterior square
of a completely continuous non-negative linear operator $A$ acting
in the ideal space $X(\Omega)$ in  terms of the spectrum of the
initial operator. The main mathematical results of this paper are
concentrated in Section~6, where we prove the existence of the
second according to the module positive eigenvalue $\lambda_2$, or
a pair of complex adjoint eigenvalues of a completely continuous
non-negative operator $A$ under the additional condition, that its
exterior square $A\wedge A$ is also nonnegative. For the case when
$A$ is a linear integral operator, the main theorem is formulated
in terms of kernels.

\section{Ideal spaces. Basic definitions and statements} Let
 $(\Omega,{\frak A},\mu)$ be a triple, consisting of
a  set $\Omega$,  $\sigma$-algebra ${\frak A}$ of measurable
subsets and  $\sigma$-finite and $\sigma$-additive complete
measure on ${\frak A}$. Denote by $S(\Omega)$ the space of all
measurable finite almost everywhere on $\Omega$ functions (further
we  consider the equivalent functions to be identical). Let
$X(\Omega)$ be a Banach ideal space, i.e. a Banach space of all
measurable on $\Omega$ functions having the following property:
from $|x_1| \leq |x_2|$, $x_1 \in S(\Omega)$, $x_2 \in X$, it
follows that $x_1 \in X$ è $\|x_1\|_X \leq \|x_2\|_X$ (the
definition and basic properties of ideal spaces are taken from
paper [10], see also [11]). Consider the support supp$X$ of the
space $X$ to be the least measurable subset, outside which all the
functions from $X$ are equal to zero. Let supp$X = \Omega $, i.e.
there exist functions from the space $X$, which are positive
almost everywhere on $\Omega$. The Banach ideal space $X$ is
called {\it regular}, if the norm in $X$ is order-continuous, i.e.
for every sequence $\{x_n\} \subset X$ from $0 \leq x_n \downarrow
0$ it follows that $\|x_n\| \rightarrow 0$. $X$ is called {\it
almost perfect}, if the norm in $X$ is order-semicontinuous, i.e.
for every sequence $\{x_n\} \subset X$ from $0 \leq x_n \uparrow x
\in X$ it follows that $\|x_n\| \rightarrow \|x\|$. It's easy to
see (see, for example, [10]), that every regular space is almost
perfect.

Let us denote by $X_{12}(\Omega \times \Omega)$ the set of all
measurable with respect to all the variables functions $x(t_1,
t_2)$ on $\Omega \times \Omega$, which satisfy the following
conditions:

(1) for almost every $t_2$ the function $x(\cdot, t_2)$ belongs to
$X$;

(2) the function $\|x(\cdot, t_2)\|$ also belongs to $X$.

If the
space $X$ is almost perfect  then the set $X_{12}$ is linear (see, for example [11], [12]). The
norm in $X_{12}$ is introduced according to the following rule:

$$\|x(t_1, t_2)\|_{12} = \|\|x(t_1, t_2)\|_{(1)}\|_{(2)},  $$
where indices (1) and (2) mean, that the norm of the space $X$ is
used firstly for the first variable, and then for the second
variable.

 The space $X_{21}(\Omega \times \Omega)$
with the norm $$\|x(t_1, t_2)\|_{21} = \|\|x(t_1,
t_2)\|_{(2)}\|_{(1)}$$ is defined similarly.

In the case of an
almost perfect $X$ both the space $X_{12}$ and the space $X_{21}$ are  almost perfect Banach ideal spaces (see
[10], [11], and also [13], p. 1234, theorem 3, where the
completeness of the space with mixed norm is proved).

 Further we will be interested in the space $\widetilde{X} = X_{12} \cap X_{21}$ of
functions that are common for the spaces $X_{12}$ and $X_{21}$. The norm in
this space is introduced by the formula: $$\|x(t_1, t_2)\|_M =
\max\{\|x(t_1, t_2)\|_{12}, \|x(t_1, t_2)\|_{21}\}.$$ Note that
the space $\widetilde{X}$ is regular if and only if the space $X$
is regular.

In connection with the introduced intersection of the spaces
$X_{12}$ and $X_{21}$ there arises a natural question of the
possibility of their coincidence. For $X = L_p$ from the Fubini
theorem it follows that $X_{12}$ and $X_{21}$ coincide according
to the fund of elements and according to their norms. However,
this is not true in  the general case. Moreover, the coincidence
of the spaces $X_{12}$ and $X_{21}$ is characteristic for the
class of the spaces $L_p$. For the regular Banach ideal space $X$
N.J. Nielsen proved, that from $X_{12} = X_{21}$ it follows that
$X$ is lattice-isomorphic to $L_p$--space (see [14]). The results
concerning this problem one can also find in [15] and [16].

\section{Tensor and exterior squares of ideal spaces} {\it The
algebraic tensor square} $X\otimes X$ of the space $X$ is defined
as the set of all functions of the form $$x(t_1,t_2) =
\sum_{i=1}^n x_{1}^{i}(t_{1}) x_{2}^{i}(t_{2}),$$ where $x_1^i,
x_2^i \in X$. Further call the elements of $X\otimes X$ {\it
degenerate functions}. By the way, {\it the algebraic exterior
square} $X\wedge X$ of the space $X$ is defined as the set of all
antisymmetric functions (i.e. functions $x(t_1,t_2)$, for which
$x(t_1, t_2) = - x(t_2, t_1)$) from $X\otimes X$.

Generally, the norm on $X\otimes X$ can be defined by different
ways. Let us go through  definitions, that  will be used
further. The norm $\alpha$ on $X\otimes X$ is called a {\it
crossnorm}, if for any $x_1, x_2 \in X$ the following equality
holds: $$ \|x_1(t_1)x_2(t_2)\|_{\alpha} =
\|x_1(t_1)\|\|x_2(t_2)\|.$$

There exists the greatest crossnorm $\pi$ (see, for example [17]),
which is defined by the equality:
 $$\|x(t_1,t_2)\|_{\pi} = \inf \sum_{i=1}^n\|x_1^i\|\|x_2^i\|,$$
where $\inf$ is taken over  all representations $
x(t_1,t_2) = \sum\limits_{i=1}^nx_1^i(t_1)x_2^i(t_2), \ x_1^i,
x_2^i \in X$.

The least crossnorm does not exist, however there exists the least
\textit{reasonable} crossnorm (the norm $\alpha$ on $X\otimes X$
is called reasonable, if $\alpha$ is a crossnorm, and the dual
norm ${\alpha}'$ is also a crossnorm on $X'\otimes X'$, where the
dual norm is defined on $X'\otimes X'$ by the equality \linebreak
$\|x'\|_{{\alpha}'} = \sup \ \{ \langle x,x' \rangle,
\|x\|_{\alpha} \leq 1\}$, $x\in X\otimes X$, $x' \in X'\otimes
X'$). The least reasonable crossnorm is denoted by $\epsilon$ and
is defined by the following rule:

$$\|x\|_{\epsilon} = \sup \ \biggl\{\biggl|\sum_{i=1}^n\langle x_1^i,x'_1
\rangle\langle x_2^i,x'_2\rangle\biggr|, \|x'_1\|_{X'},\|x'_2\|_{X'}
\leq 1\biggr\}.$$

Note that the completion of the algebraic tensor square $X\otimes
X$ of the ideal space $X$ with respect  to  the norms $\pi$ or
$\epsilon$ will not be ideal. It is natural to define the norm on
$X\otimes X$ in such a way that the completion with respect to this norm will  be an ideal space or
it's part. With this purpose V.L. Levin introduced in [18] the
following crossnorm on $X\otimes X$:

$$\|x\|_{L} = \inf \biggl\{\|u\|_X: u \geq \biggl|\sum_{i} x_{1}^{i} \langle
x_{2}^{i}, x^* \rangle\biggr|, \|x^*\|_{X^*} \leq 1 \biggr\}.$$

 V.L. Levin proved (see [18], p. 55, proof of the theorem 1),
that the topology of the space $X_{21}$ induces on $X\otimes X$
the same topology, as  the norm $L$.

Note that every norm $\alpha$ on $X\otimes X$ which satisfies
inequalities $ \epsilon \leq \alpha \leq \pi$ will be reasonable.
In particular, the Levin's norm $L$ is  reasonable (see [18],
p. 53, lemma 3).

Further let us call the completion of $X\otimes X$ with respect  to the
Levin's norm the {\it $L$-complete} tensor square of the space $X$
 and let us denote it by $\widetilde{(X\otimes X)}_{L}$. As it was
noticed above, $\widetilde{(X\otimes X)}_{L}$ is a closed subspace
of $X_{21}$. The space $\widetilde{(X\otimes X)}_{L}$ was studied
thoroughly by A.V. Bukhvalov in papers [19] and [13]. He also
proved the criteria of the coincidence of this space and $X_{21}$
(see [13], p. 7, theorem 0.1, and also [19], p. 1235, theorem 4).

{\bf Bukhvalov's theorem.} \textit{Let $X$ be a Banach ideal space.
Then the following statements are equivalent:
\begin{enumerate}\parskip=-1mm
 \item[\rm (i)] the set of all degenerate functions is dense in $X_{21}$;
 \item[\rm (ii)] the equality $\widetilde{(X \otimes X)}_{L} =
X_{21}$ is true;
 \item[\rm (iii)] the space $X$ is regular.\end{enumerate}}

 Bukhvalov's theorem implies that  for a regular Banach
ideal space the following equality holds:
$$
\widetilde{(X\otimes
X)}_{L}^{a} = X_{21}^a(\Omega \times \Omega),
$$
where
$X_{21}^a(\Omega \times \Omega)$ is a subspace of the space
$X_{21}(\Omega \times \Omega)$, which consists of antisymmetric
functions.

 Further note that $\wedge$-product
of arbitrary functions $x_1, x_2 \in X \ $ $(x_1 \wedge
x_2)(t_1,t_2) = x_1(t_1)x_2(t_2) - x_1(t_2)x_2(t_1)$ belongs to
the space $X_{21}(\Omega \times \Omega)$, and the following
equality holds:
$$
x_1 \wedge x_2(t_1,t_2) = - x_1 \wedge
x_2(t_2,t_1).
$$
Therefore, for any $x_1, x_2 \in X $ the
function $x_1 \wedge x_2$ belongs to the subspace $X_{21}^a(\Omega
\times \Omega)$.

Moreover, {\it an arbitrary antisymmetric function $x(t_1,t_2)$,
which belongs to the space $X_{21}$, at the same time belongs to
$X_{12}$, and the following equality holds} $$\|x(t_1,t_2)\|_{21}
= \|x(t_1,t_2)\|_{12}.$$

 Really, let $x(t_1, t_2)$ be equal to $ - x(t_2,
t_1)$. In this case $$\|x(t_1,t_2)\|_{21} =
\|\|x(t_1,t_2)\|_{(2)}\|_{(1)} = \|\|- x(t_2,t_1)\|_{(2)}\|_{(1)}
= $$ $$ = \|\|x(t_1,t_2)\|_{(1)}\|_{(2)} = \|x(t_1,t_2)\|_{12}.$$
That is why further we will assume that  the completion of the algebraic
exterior square $X\wedge X$ of the space $X$ is taken with respect to the
symmetric crossnorm
$$
\|x(t_1, t_2)\|_M = \max\{\|x(t_1,
t_2)\|_{12}, \|x(t_1, t_2)\|_{21}\}.
$$
This completion
coincides with the closed subspace of antisymmetric functions of
the space \linebreak $\widetilde{X}(\Omega \times \Omega)$.

The subspace $\widetilde{X}^a(\Omega \times \Omega)$ is isomorphic
in the category of Banach spaces to the space $\widetilde{X}(W)$,
where $W$ is the measurable subset $\Omega \times \Omega$, for
which the sets $W \cap \widetilde W$ and $(\Omega \times \Omega)
\setminus (W \cup \widetilde W)$ have zero measure; here
$\widetilde W = \{(t_2,t_1): \ (t_1,t_2) \in W\}$ (such sets do
always exist). Really, extending the functions from
$\widetilde{X}(W)$ as antisymmetric functions to $\Omega \times
\Omega$, we obtain the set of all the functions from
$\widetilde{X}^a(\Omega \times \Omega)$. Further, setting the norm
of a function in $\widetilde{X}(W)$ to be equal to the norm of its
extension, we get that the spaces $\widetilde{X}^a(\Omega \times
\Omega)$ and $\widetilde{X}(W)$ are isomorphic in the category of
normed spaces.

Therefore, for an almost perfect Banach ideal space the following
equality holds:
$$
\widetilde{(X\wedge X)}_{M} =
\widetilde{(X\otimes X)}_{M}^{a} = \widetilde{X}_d(W),
$$
where
$\widetilde{X}_d(W)$ is the closure of the set of all degenerate
functions from $\widetilde{X}(W)$ in the norm of
$\widetilde{X}(W)$.

\section{Tensor and exterior squares of linear operators in ideal
spaces} Let $A, B$ be continuous linear operators acting in the
ideal space $X$. Define the {\it algebraic tensor product} of the
operators $A$ and $B$ as the operator  $A \otimes B$ in the space
$X\otimes X$, defined on degenerate functions by the equality $$(A
\otimes B)x(t_1,t_2) = \sum_j Ax_1^{j}(t_1) \cdot Bx_2^{j}(t_2)
\qquad \bigg(x(t_1,t_2) = \sum_j x_1^{j}(t_1) \cdot
x_2^{j}(t_2)\bigg).$$

 A crossnorm $\alpha$ on $X\otimes X$ is called
\textit{quasiuniform}, if for any continuous linear operators
$A,B$ one has:
$$
\|A \otimes B\|_{\alpha} \leq c\|A\|\|B\|,
$$
where $c$ is some constant. If $c = 1$, then such  crossnorm is
called \textit{uniform}.

The greatest and the least reasonable norms $\pi$ and $\epsilon$
are uniform. However  Levin's norm $L$ in general case will not
be quasiuniform. That is why we need the following definition.

Define the {\it tensor product} of the operators $A$ and $B$ (if
it exists) as the linear operator $A \otimes B$ in the space
$\widetilde{(X\otimes X)}_{L} $, defined on degenerate functions
by the equality $$(A \otimes B)x(t_1,t_2) = \sum_j Ax_1^{j}(t_1)
\cdot Bx_2^{j}(t_2) \qquad \bigg(x(t_1,t_2) = \sum_j x_1^{j}(t_1)
\cdot x_2^{j}(t_2)\bigg),$$ and on arbitrary functions by
extension via continuity from the subspace of degenerate functions
onto the whole of $\widetilde{(X\otimes X)}_{L}$ (if the
extension is  bounded).

Let us formulate the following statement concerning with the
tensor product of two operators in ideal spaces (see [18], p. 62,
proposition 6).

{\bf Levin's theorem.} \textit{Let $X$ be a Banach ideal space,
$A: X \rightarrow X$ be a regular operator, i.e. an operator,
which can be represented in the form $A = A_1 - A_2$, where $A_1$
and $A_2$ are nonnegative linear operators (with respect to the cone of nonnegative
functions), and let  $B: X \rightarrow  X$ be a continuous
linear operator. Then the
tensor product $A \otimes B$ does
exist.}

Further, when studying the spectral properties of the tensor square $A
\otimes A$ of the operator $A$, we will have to impose conditions on the
operator $A$,
that  are stronger than nonnegativity or regularity. Let us give the following definition. A linear
operator $A:X \rightarrow X$ is called {\it resolvent-regular}, if for any
$\lambda$, which is not in $\sigma(A)$, the resolvent operator
$R(\lambda, A) = (\lambda I - A)^{-1}$ is regular. The class of
resolvent-regular operators includes, for example, Hilbert-Schmidt
operators and operators, such that their certain  power  is a Hilbert-Schmidt
operator.

 Further let us examine the operator $A \wedge
A$, defined as the restriction of the operator $A \otimes A$ onto
the subspace $\widetilde{(X\otimes X)}_{M}^{a}$. It is obvious
that for degenerate antisymmetric functions the operator $A \wedge
A$ can be defined by the equality $$(A \wedge A)x(t_1,t_2) =
\sum_j Ax_1^{j}(t_1) \wedge Ax_2^{j}(t_2) \qquad x(t_1,t_2) =
\sum_j x_1^{j}(t_1) \wedge x_2^{j}(t_2).$$

\section{ Spectrum of the tensor square of linear operators in ideal
spaces} As usual, we denote by $\sigma(A)$ the spectrum of an
operator $A$, and we denote by $\sigma_{p}(A)$ the point spectrum, that is,
the set of all eigenvalues of the operator $A$. Denote by
$\sigma_{eb}(A)$ the Browder essential spectrum of the operator
$A$. Thus $\sigma(A)\setminus \sigma_{eb}(A)$ will be the set of
all isolated finite dimensional eigenvalues of the operator $A$
(for more detailed information see [20], [21]).

In papers [20]--[23] T. Ichinose obtained the results,
representing the spectra and the parts of the spectra of the
tensor products of linear bounded operators in terms of the
spectra and parts of the spectra of the  operators given under the
natural assumptions, that the corresponding crossnorm is
reasonable and quasiuniform. Among the mentioned results there are
the explicit formulae, expressing the set of all isolated finite
dimensional eigenvalues and the Browder essential spectrum of the
operator $A \otimes A$ in terms of the parts of the spectrum of
the  operator given (see [20], p. 110, theorem 4.2). In particular,
Ichinose proved, that for the tensor square of a linear bounded
operator $A$ the following equalities hold: $$\sigma(A\otimes A) =
\sigma(A) \, \sigma(A);$$
 $$\sigma(A\otimes A) \setminus
\sigma_{eb}(A\otimes A) = (\sigma(A) \setminus \sigma_{eb}(A)) \,
(\sigma(A) \setminus \sigma_{eb}(A)) \setminus
(\sigma_{eb}(A)\sigma(A));\eqno(1)$$ $$\sigma_{eb}(A\otimes A) =
\sigma_{eb}(A) \, \sigma(A).$$

 For a completely continuous operator the following
equalities hold: $$(\sigma(A)\setminus \sigma_{eb}(A))\setminus
\{0\} = \sigma_{p}(A)\setminus \{0\}; \ \sigma_{eb}(A) = \{0\} \
\mbox{or} \ \emptyset. $$ So from (1) we can get the complete
information about the nonzero eigenvalues of the tensor square of
a completely continuous operator: $$\sigma_{p}(A\otimes
A)\setminus \{0\} = (\sigma_{p}(A) \, \sigma_{p}(A))\setminus
\{0\}. \eqno(2)$$ Here zero can be either a finite or infinite
dimensional eigenvalue of $A\otimes A$ or a point of the essential
spectrum.

In paper [9] there have been examined the case, when a  linear
operator acts in the space $L_p(\Omega)$ ($C(\Omega)$), and the
corresponding crossnorm is reasonable and quasiuniform. That is
why the formula for the spectrum of $A\otimes A$ directly follows
from the results of T. Ichinose.

However in general case the crossnorm $L$ is reasonable, but not
quasiuniform (see [24]), and therfore we need a different  proof for the
statement about the spectrum of $A\otimes A$. The proof, given
below, is based on the reasoning, made by A.S. Kalitvin (see [11],
p. 83, theorem 3.10), for the case of the operator $A\otimes I +
I\otimes A$ in a regular ideal space.

 {\bf Theorem 1.} \textit{Let $X$
be an almost perfect Banach ideal space, and let $A:X \rightarrow X$
be a completely continuous nonnegative with respect to the cone of
nonnegative functions in $X$ resolvent-regular operator. Then for
the point spectrum of the operator $A\otimes A$, acting in the
space $\widetilde{(X\otimes X)}_{L}$, the following equality
holds: $$\sigma_{p}(A\otimes A)\setminus \{0\} = (\sigma_{p}(A) \,
\sigma_{p}(A))\setminus \{0\}.$$ }

\vspace{-0.3cm}

{\bf Proof.} Let us examine the operators $A \otimes I$ and $I
\otimes A$, acting in $\widetilde{(X\otimes X)}_{L}$. Let us
prove, that the following inclusions are true: $\sigma(A \otimes I)
\subseteq \sigma(A)$ and $\sigma(I \otimes A) \subseteq \sigma(A)$.
Let us prove the first inclusion (the second inclusion can be proved
by analogy). Let $\lambda$ does not belong to $\sigma(A)$.
Then the operator $\lambda I - A$ is invertible. Let us define the
operator $(\lambda I - A)^{-1}\otimes I$ on $X \otimes X$. Since
the operator $(\lambda I - A)^{-1}$ is regular, we can apply the
Levin's theorem. As it follows from  Levin's theorem, the
operator $(\lambda I - A)^{-1}\otimes I$ can be extended from $X
\otimes X$ onto the whole of $\widetilde{(X\otimes X)}_{L}$. It is
easy to see, that the operator $(\lambda I - A)^{-1}\otimes I$ is
inverse on $X \otimes X$ for the operator $\lambda I - (A \otimes
I)$. So, its extension $(\lambda I - A)^{-1} \otimes I$
will be inverse for the operator $\lambda I - (A \otimes I)$ on
the whole of $\widetilde{(X\otimes X)}_{L}$. That is why $\lambda$
does not belong to $\sigma(A\otimes I)$, and the inclusion
$\sigma(A \otimes I) \subseteq \sigma(A)$ is proved.

Thus  as $A\otimes A = (A \otimes I)(I \otimes A)$ and the operators
$A \otimes I$ and $I \otimes A$ are, obviously, commuting,
the following relation is true
$$
\sigma(A\otimes A) = \sigma((A
\otimes I)(I \otimes A)) \subseteq \sigma(A \otimes I)\sigma(I
\otimes A).
$$
Now,  applying the inclusions $\sigma(A \otimes I)
\subseteq \sigma(A)$ and $\sigma(I \otimes A) \subseteq
\sigma(A)$, proved above, we see that the following inclusion is  true as well
$$
\sigma(A\otimes A) \subseteq \sigma(A)\sigma(A).
$$
 since due to the
complete continuity of the operator $A$, its spectrum, except,
probably, zero, consists of isolated finite dimensional
eigenvalues. That is why the following relations hold:
$$\sigma_p(A\otimes A)\setminus \{0\} \subseteq \sigma(A\otimes
A)\setminus \{0\} \subseteq (\sigma(A)\sigma(A))\setminus \{0\}
=$$ $$= (\sigma_p(A)\sigma_p(A))\setminus \{0\},$$ i.e. we proved:
$$\sigma_p(A\otimes A)\setminus \{0\} =
(\sigma_p(A)\sigma_p(A))\setminus \{0\}.$$

 Now, let us prove the reverse inclusion. For this we will examine the
extension $\widetilde {(A \otimes A)}_\epsilon$ of the operator $A
\otimes A$ onto the whole $\widetilde{(X \otimes X)}_\epsilon$,
where $\widetilde{(X \otimes X)}_\epsilon$ is a completion $X
\otimes X$ with respect to  the "weak" \ crossnorm $\epsilon$. As it follows
from the results of J.R. Holub (see [25], p. 401, theorem 2), the
operator $\widetilde {(A \otimes A)}_\epsilon$ is completely
continuous in $\widetilde{(X \otimes X)}_\epsilon$. Let us prove, that $
\sigma_{p}(\widetilde {(A \otimes A)}_\epsilon)\setminus \{0\}
\subseteq \sigma_p(A\otimes A)\setminus \{0\}.$
To check  this it is enough  to prove, that any eigenfunction of the operator $\widetilde
{(A \otimes A)}_\epsilon$, corresponding to a nonzero eigenvalue,
belongs to the space $\widetilde{(X\otimes X)}_{L}$. Let $\lambda$
be an arbitrary nonzero eigenvalue of the operator $\widetilde {(A
\otimes A)}_\epsilon$. Formula (3) implies that there
exist indices $i, \ j$, for which $\lambda = \lambda_i \lambda_j$
(here $\{\lambda_k\}$ is the set of all nonzero eigenvalues of the
operator $A$, enumerated without regard to multiplicity). Since
$\lambda$ is an isolated finite dimensional eigenvalue, $i, \ j$
can take  only finite number of different values. Let us
enumerate all the pairs of such values. Let $\lambda =
\lambda_{i_k} \lambda_{j_k} \ \ (k = 1, \ldots, p)$. Decompose the
space $X$ into the direct sum of subspaces: $$X = X_1 \oplus
\ldots \oplus X_p \oplus R,$$ where $X_k = \ker (A -
\lambda_{i_k})^{m_k}$, $m_k$ are the multiplicities of
$\lambda_{i_k}$. Under this decomposition  $\widetilde{(X \otimes
X)}_\epsilon$  also decomposes into the direct sum of subspaces:
$$\widetilde{(X \otimes X)}_\epsilon = (X_1 \otimes X_1) \oplus
\ldots \oplus (X_1 \otimes X_p)\oplus \ldots \oplus(X_p \otimes
X_1) \oplus \ldots \oplus (X_p \otimes X_p)\oplus $$ $$ \oplus
(X_1 \otimes R) \oplus \ldots \oplus (X_p \otimes R)\oplus(R
\otimes X_1) \oplus \ldots \oplus (R \otimes X_p)\oplus
\widetilde{(R \otimes R)}_\epsilon.$$ Since $X_l \otimes X_m$,
$X_l \otimes R$, $R \otimes X_l$, $\widetilde{(R \otimes
R)}_\epsilon$ are  invariant subspaces for the operator
$\widetilde{(A\otimes A)}_\epsilon$, the following equality holds:
$$\sigma \widetilde{(A\otimes A)}_\epsilon = $$ $$ = \bigcup_{l,m}
(\sigma(A\otimes A, X_l \otimes X_m)\cup \sigma(A\otimes A, R
\otimes X_m)\cup \sigma(A\otimes A, X_l \otimes R))\cup
\sigma(A\otimes A,\widetilde{R \otimes R}_\epsilon),$$ where the
notation $\sigma(A\otimes A, X_l \otimes X_m)$ means the spectrum
of the restriction of the operator $A\otimes A$ onto the
corresponding subspace. Since $X_l \otimes X_m$, $R \otimes X_m$,
$X_l \otimes R$, $\widetilde{(R \otimes R)}_\epsilon$ are the spaces
with uniform crossnorms, we can apply the results of T. Ichinose,
and therefore the following equalities hold: $$\sigma(A\otimes A, X_l
\otimes X_m) = \sigma(A, X_l)\sigma(A,X_m); \ \sigma(A\otimes A,
X_l \otimes R) = \sigma(A, X_l)\sigma(A,R);$$ $$\sigma(A\otimes A,
R \otimes X_m) = \sigma(A, R)\sigma(A,X_m); \ \sigma(A\otimes A,
\widetilde{(R \otimes R)}_\epsilon) = \sigma(A, R)\sigma(A, R). $$
Since $\lambda_{i_k}$ and $\lambda_{j_k}$ do not belong $\sigma(A,
R)$ for any values of indices $k \ \ (k = 1, \ldots, p)$,
 $\lambda$ does not belong to $ \bigcup_{l,m} (\sigma(A\otimes A, R
\otimes X_m)\cup \sigma(A\otimes A, X_l \otimes R))\cup
\sigma(A\otimes A,\widetilde{(R \otimes R)}_\epsilon).$ As it
follows, $\lambda \in  \bigcup_{l,m} (\sigma(A\otimes A, X_l
\otimes X_m))$. Further it is obvious, that for any $l,m \ \ ( 1
\leq l,m \leq p) \ $ $X_l \otimes X_m$ belongs to the algebraic
tensor square $X \otimes X$ and therefore it  belongs to the
space $\widetilde{(X\otimes X)}_{L}$. So,  for an arbitrary $
\lambda \in \sigma_{p}(\widetilde {(A \otimes
A)}_\epsilon)\setminus \{0\}$ the inclusion $\lambda \in
\sigma_p(A\otimes A)\setminus \{0\}$ is true. $\square$

Let us notice, that under conditions of  Theorem 1 the
inclusion:
 $$\sigma_{eb}(A \otimes A) \subseteq \sigma_{eb}(A)\sigma_{eb}(A) = \{0\}$$
follows from its proof. Moreover, for an arbitrary $\lambda \in
(\sigma_{p}(A\otimes A)\setminus \{0\})$ the following equality
holds:
 $$\ker (A\otimes A - \lambda I\otimes I) = \bigoplus_{k}\ker (A - \lambda_{i_k} I)\otimes \ker (A - \lambda_{j_k}I),$$
where the summation is taken over all the numbers $k$ of the pairs
$\lambda_{i_k},\lambda_{j_k}$, for which
$\lambda_{i_k},\lambda_{j_k} \in (\sigma_{p}(A)\setminus \{0\})$
and $\lambda = \lambda_{i_k}\lambda_{j_k}$.

As it follows from the results of the paper [26], the statements
of  Theorem 1 will be also true  for the case, when the tensor
square $A\otimes A$ of the operator $A$ acts in the space
$\widetilde{(X\otimes X)}_{M}$.

\section{Spectrum of the exterior square of linear operators in
ideal spaces} For the exterior square, which is the restriction of
the tensor square, the following inclusions hold:
$$\sigma(A\wedge A) \subset \sigma(A\otimes A)$$
 $$\sigma_{p}(A\wedge A) \subset
\sigma_{p}(A\otimes A).$$

Applying Theorem 1, we can prove the following statement:

{\bf Theorem 2.} \textit{Let $X$ be an almost perfect ideal space
and let $A: X \rightarrow X$ be a completely continuous nonnegative
with respect to the cone of nonnegative functions in $X$
resolvent-regular operator. Let $\{\lambda_{i}\}$ be the set of
all nonzero eigenvalues of the operator $A$, repeated according to
multiplicity. Then all the possible products of the type
$\{\lambda_{i} \lambda_{j} \}$, where $i < j$, form the set of all
the possible (except, probably, zero) eigenvalues of the exterior
square of the operator $A \wedge A $, repeated according to
multiplicity}.

The proof copies the proof of the corresponding statement from [9]
(see [9], p. 12, Theorem 1).

\section{ Generalization of the Gantmakher--Krein theorems in the
case of $2$-totally nonnegative operators in ideal spaces} Let us
prove some generalizations of the Gantmakher--Krein theorems in
the case of operators in ideal spaces, using the generalized
Krein--Rutman theorem (see [27]) about completely continuous
operators, leaving invariant an almost reproducing cone $K$ in a
Banach space (for such operators the spectral radius $\rho(A)$
belongs to $\sigma_{p}(A)$).

{\bf Theorem 3.} \textit{Let $X(\Omega)$ be an almost perfect
ideal space, and, respectively, let $\widetilde{X}_d(W)$ be the
closure of the set of all degenerate functions in intersection of
the spaces with mixed norms. Let $A: \ X(\Omega) \rightarrow
X(\Omega)$ be a completely continuous nonnegative with respect to
the cone of nonnegative functions in $X$ resolvent-regular
operator, such that  $\rho(A)
> 0$. Let the exterior square $A \wedge A: \ \widetilde{X}_d(W)
\rightarrow \widetilde{X}_d(W)$ be nonnegative with respect to the
cone of nonnegative functions in $\widetilde{X}_d(W)$, and
$\rho(A \wedge A)>0$. Then the operator $A$ has a positive
eigenvalue $\lambda_{1} = \rho(A)$. Moreover, if there is only one
eigenvalue on the spectral circle $|\lambda| = \rho(A)$, then the
operator $A$ has the second positive eigenvalue $\lambda_2 <
\lambda_{1}$. If there is more than one eigenvalue on the spectral
circle $|\lambda| = \rho(A)$, then either there is at least one
pair of complex adjoint among them, or $\lambda_{1}$ is a multiple
eigenvalue}.

 {\bf Proof.} Enumerate eigenvalues of a completely continuous operator $A$,
repeated according to multiplicity, in order of decrease of their
modules:
$$|\lambda_{1}| \geq | \lambda_{2}| \geq |\lambda_{3}| \geq  \ldots$$
Applying  the generalized Krein--Rutman theorem to $A$ we get:
$\lambda_{1} = \rho(A)>0$. Now applying the generalized
Krein--Rutman theorem to the operator $A \wedge A$ we get:
$\rho(A\wedge A) \in \sigma\!_{p}(A\wedge A)$. (Note, that from
$\sigma_{eb}(\wedge^2 A)
\subseteq\sigma_{eb}(\widetilde{(\otimes^2 A)}_{M_2}) \subseteq
\sigma_{eb}(A)\sigma_{eb}(A) = \{0\}$, it follows, that the
operator $A\wedge A$ also satisfy the conditions of the
generalized Krein--Rutman theorem).

 As it follows from the statement of  theorem 2, the exterior
 square of the operator
$A$  has no other nonzero eigenvalues, except for all the possible
products of the form $\lambda_{i}\lambda_{j}$, where $ i < j$. So,
we  conclude  that $\rho(A\wedge A)>0$ can be represented
in the form of the product $\lambda_{i}\lambda_{j}$ with some
values of the indices $i,j$, \ $i < j$. Thus,  if there is only one
eigenvalue on the spectral circle $|\lambda| = \rho(A)$, from the
fact that eigenvalues are numbered in decreasing order it follows
that
 $\rho(A\wedge A) =  \lambda_{1}\lambda_{2}$.
  Therefore  $\lambda_{2} = \frac{\rho(A\wedge
A)}{\lambda_{1}}>0$.

If there is $m \ (m \geq 2)$ eigenvalues on the spectral circle
$|\lambda| = \rho(A)$, then $\rho(A\wedge A) =
\lambda_{i}\lambda_{j},$ where $1\leq i < j \leq m$. So, both
$\lambda_{i}$ and $\lambda_{j}$ are situated on the spectral
circle $|\lambda| = \rho(A)$, and from the positivity of their
product it follows, that $\lambda_{i}$ è $\lambda_{j}$ are either
a pair of complex adjoint eigenvalues, or both are positive and
coincide with $\rho(A)$. $\square$

It is well-known (see, for example, [25], p. 55, corollary from
the proposition 2.1), that a linear integral operator $A$, acting
in the Banach ideal space $X(\Omega)$, is nonnegative if and only
if its kernel $k(t,s)$ is nonnegative almost everywhere on
$\Omega$. It is also well-known (see [9]), that the exterior power
of a linear integral operator can be considered as a linear
integral operator, acting in the space $\widetilde{X}_d(W)$ with
the kernel equal to the second associated to the kernel of the
 operator given. That is why it is not difficult to reformulate
Theorem 3 in terms of kernels of linear integral operators. In
this case the conditions of Theorem 3 can  be
easily verified.

{\bf Theorem 4.} \textit{Let a completely continuous
resolvent-regular linear integral operator $A$ act in an almost
perfect ideal space $X(\Omega)$. Let the kernel $k(t,s)$ of the
operator $A$ be nonnegative and not equal identically to zero
almost everywhere on the Cartesian square $\Omega \times \Omega$.
Let the second associated kernel $k\wedge k(t_1, \ t_2, \ s_1, \
s_2)$ be nonnegative and not equal identically to zero almost
everywhere on the Cartesian square $W \times W$, where $W$ is a
measurable subset, possessing the following properties:}

1) $\mu (W \cap \widetilde W) = 0$;

2) $\mu((\Omega \times \Omega) \setminus (W \cup \widetilde W)) =
0. $ \ \ \ \ \ \ \ \ \ $(\widetilde W = \{(t_2,t_1): \ (t_1,t_2)
\in W\})$

\textit{ Then the operator $A$ has a positive eigenvalue
$\lambda_{1} = \rho(A)$. Moreover, if there is only one eigenvalue
on the spectral circle $|\lambda| = \rho(A)$, then the operator
$A$ has the second positive eigenvalue $\lambda_2 < \lambda_{1}$.
If there is more than one eigenvalue on the spectral circle
$|\lambda| = \rho(A)$, then either there is at least one pair of
complex adjoint among them, or $\lambda_{1}$ is a multiple
eigenvalue}.

Note that in Theorem 4 the kernel is not presupposed to be
continuous, we   assume only, that the operator $A$ acts in one of
almost perfect ideal spaces.

Moreover, Theorem 3 can be generalized in the case, when the
exterior square $A \wedge A$ of the operator $A$ leaves invariant
an arbitrary almost reproducing cone $\widetilde{K}$ in
$\widetilde{X}_d(W)$. But in this case certain difficulties,
related to the  testing of the assumption of the generalized
theorem, can arise.

\newpage

\section*{ References}

\medskip

1. {\it Gantmacher F.R., Krein M.G.} Oscillation Matrices and
Kernels and Small Vibrations of Mechanical Systems.
--- AMS Bookstore, 2002. --- 310 p.

2. {\it Kellog O.D.} Orthogonal function sets arising from
integral equations // American Journal of Mathematics. -- 1918. --
Vol. 40. -- P. 145-154.

3. {\it Pinkus A.} Spectral properties of totally positive kernels
and matrices // Total positivity and its applications / M. Gasca,
C.A. Micchelli. -- Dordrecht, Boston, London: Kluwer Acad. Publ.,
1996. -- P. 1-35.

4. {\it Karlin S.} Total positivity. -- Stanford University Press,
California, 1968. -- Vol. 1. -- 576 p.

5. {\it Eveson S.P.} Eigenvalues of totally positive integral
operators // Bull. London Math. Soc. -- 1997. -- Vol. 29. -- P.
216-222.

6. {\it Sobolev A.V.} Abstract oscillatory operators
// Proceedings of the Seminar on differential equations. -- Kuibishev, 1977.
-- No. 3. -- P. 72-78 (Russian).

7. {\it Yudovich V.I.} Spectral properties of an evolution
operator of a parabolic equation with one space variable and its
finite-dimensional analogues // Uspekhi Mat. Nauk. -- 1977. --
Vol.32, No. 1. -- P. 230-232 (Russian).

8. {\it Pokornyi Yu.V., Penkin O.M., Pryadiev V.L., Borovskikh
A.V., Lazarev K.P., Shabrov S.A.} Differential equations on
geometrical graphs.
--- Moscow: FIZMATLIT, 2004. --- 272 p.

9. {\it Kushel O.Y., Zabreiko P.P.} Gantmakher-Krein theorem for
2-nonnegative operators in spaces of  functions // Abstract and
Applied Analysis. -- 2006. -- Article ID 48132. -- P. 1-15.

10. {\it Zabreiko P.P.} Ideal spaces of functions // Vestnik
Yaroslav. Univ. -- 1974. -- No. 4. -- P. 12-52 (Russian).

11. {\it Kalitvin A.S.} Linear operators with partial integrals.
--- Voronezh, 2000. --- 252 p. (Russian).

12. {\it Kantorovich L.V., Akilov G.P.} Functional Analysis. ---
2nd rev. Moscow, 1977 (Russian). English transl.: Pergamon Press,
Oxford, 1982.

13. {\it Bukhvalov A.V.} On spaces with mixed norm // Vestnik
Leningrad. Univ. -- 1973. -- No. 19. -- P. 5–12 (Russian); English
transl.: Vestnik Leningrad Univ. Math. -- 1979. -- No. 6. -- P.
303–311.

14. {\it Nielsen N.J.} On Banach ideals determined by Banach
lattices and their applications // Diss. Math. -- 1973. ¹ 109.

15. {\it Bukhvalov A.V.} Generalization of the Kolmogorov-Nagumo
theorem on the tensor product // Qualitative and Approximate
Methods for the Investigation of Operator Equations. --
Jaroslavl', 1979. -- No. 4. -- P. 48–65 (Russian).

16. {\it Boccuto A., Bukhvalov A.V., Sambucini A.R.} Some
inequalities in classical spaces with mixed norms // Positivity.
-- 2002. -- Vol. 6. No. 4. -- P. 393-411.

17. {\it Tsoy-Wo Ma.} Classical analysis on normed spaces. ---
World Scientific Publishing, 1995.

18. {\it Levin V.L.} Tensor products and functors in categories of
Banach spaces defined by KB-lineals // Trudy Moskov. Mat. Obshch.
-- 1969. -- Vol. 20. -- P. 43–81 (Russian); English transl.:
Trans. Moscow Math. Soc. -- 1969. -- Vol. 20. -- P. 41–78.

19. {\it Bukhvalov A.V.} Vector-valued function spaces and tensor
products // Siberian Math. J. -- 1972. -- Vol. 13, No. 6. -- P.
1229-1238 (Russian).

20. {\it Ichinose T.} Spectral properties of tensor products of
linear operators. I // Transactions of the American Mathematical
Society. -- 1978. -- Vol. 235. -- P. 75-113.

21. {\it Ichinose T.} Spectral properties of tensor products of
linear operators. II // Transactions of the American Mathematical
Society. -- 1978. -- Vol. 237. -- P. 223-254.

22. {\it Ichinose T.} Operators on tensor products of Banach
spaces
// Trans. Am. Math. Soc. -- 1972. -- Vol. 170. -- P. 197-219.

23. {\it Ichinose T.} Operational calculus for tensor products of
linear operators in Banach spaces. // Hokkaido Math. J. -- 1975.
-- Vol. 4. No. 2. -- P. 306-334.

24. {\it Bukhvalov A.V.} Application of methods of the theory of
order-bounded operators to the theory of operators in $L^p$-spaces
// Uspekhi Mat. Nauk. -- 1983. --  Vol. 38, Issue 6(234). -- P. 37–83.

25. {\it Holub J.R.} Compactness in topological tensor products
and operator
 spaces // Proc. Am. Math. Soc. -- 1972. -- Vol. 36. -- P. 398-406.

26. {\it Zabreiko P.P.} The spectrum of linear operators acting in
various Banach spaces
// Qualitative and Approximate Methods for the
Investigation of Operator Equations. -- Jaroslavl', 1976. -- No.
1. -- P. 39-47 (Russian).

27. {\it Zabreiko P.P. and  Smitskikh S.V.} A theorem of M. G.
Krein and M. A. Rutman // Functional Analysis and Its
Applications. -- 1980. -- Vol. 13. -- P. 222–223.

\end{document}